\documentclass[oneside,american]{amsart}
\usepackage[T1]{fontenc}
\usepackage[utf8x]{inputenc}
\usepackage{babel}
\usepackage{prettyref}
\usepackage{mathtools}
\usepackage{amsthm}
\usepackage{amstext}
\usepackage{graphicx}
\usepackage[unicode=true,pdfusetitle,
 bookmarks=true,bookmarksnumbered=false,bookmarksopen=false,
 breaklinks=false,pdfborder={0 0 1},backref=false,colorlinks=false]
 {hyperref}

\makeatletter

\numberwithin{equation}{section}
\numberwithin{figure}{section}
 \theoremstyle{definition}
 \newtheorem*{defn*}{\protect\definitionname}

\usepackage{amsfonts}
\usepackage{dsfont,mathtools}
\usepackage[format=hang]{caption}
\usepackage{prettyref}

\newrefformat{prop}{Proposition \ref{#1}}
\newrefformat{lem}{Lemma \ref{#1}}
\newrefformat{thm}{Theorem \ref{#1}}
\newrefformat{cor}{Corollary \ref{#1}}
\newrefformat{rem}{Remark \ref{#1}}
\newrefformat{example}{Example \ref{#1}}
\newrefformat{def}{Definition \ref{#1}}
\newrefformat{eq}{(\ref{#1})}

\makeatother

  \providecommand{\definitionname}{Definition}

\begin{document}

\title[Attraction-Based Computation of Hyperbolic LCS]{Attraction-Based Computation of Hyperbolic Lagrangian Coherent Structures}

\author{Daniel Karrasch}

\address{ETH Zürich, Institute of Mechanical Systems, Leonhardstrasse 21, 8092
Zurich, Switzerland}

\email{kadaniel@ethz.ch}

\author{Mohammad Farazmand}

\address{ETH Zürich, Institute of Mechanical Systems, Leonhardstrasse 21, 8092
Zurich, Switzerland}

\address{ETH Zürich, Department of Mathematics, Rämistrasse 101, 8092 Zurich,
Switzerland}

\email{farazmam@ethz.ch}

\author{George Haller}

\address{ETH Zürich, Institute of Mechanical Systems, Leonhardstrasse 21, 8092
Zurich, Switzerland}

\email{georgehaller@ethz.ch}

\date{\today}



\begin{abstract}
Recent advances enable the simultaneous computation of both attracting
and repelling families of Lagrangian Coherent Structures (LCS) at
the same initial or final time of interest. Obtaining LCS positions
at intermediate times, however, has been problematic, because either
the repelling or the attracting family is unstable with respect to
numerical advection in a given time direction. Here we develop a new
approach to compute arbitrary positions of hyperbolic LCS in a numerically
robust fashion. Our approach only involves the advection of attracting
material surfaces, thereby providing accurate LCS tracking at low
computational cost. We illustrate the advantages of this approach
on a simple model and on a turbulent velocity data set.
\end{abstract}
\maketitle

\section{Introduction}

Hyperbolic Lagrangian Coherent Structures (LCS) in a two-dimensional
unsteady flow are locally most repelling or most attracting material
lines over a given finite time interval $I=[t_{1},t_{2}]$ of interest
\cite{Haller2000}. Both mathematical methods and intuitive diagnostic
tools have been developed to locate LCS in finite-time unsteady velocity
fields with general time dependence (see \cite{Haller2015} for a
recent review.)

\begin{figure}
\centering
\includegraphics[width=0.45\textwidth]{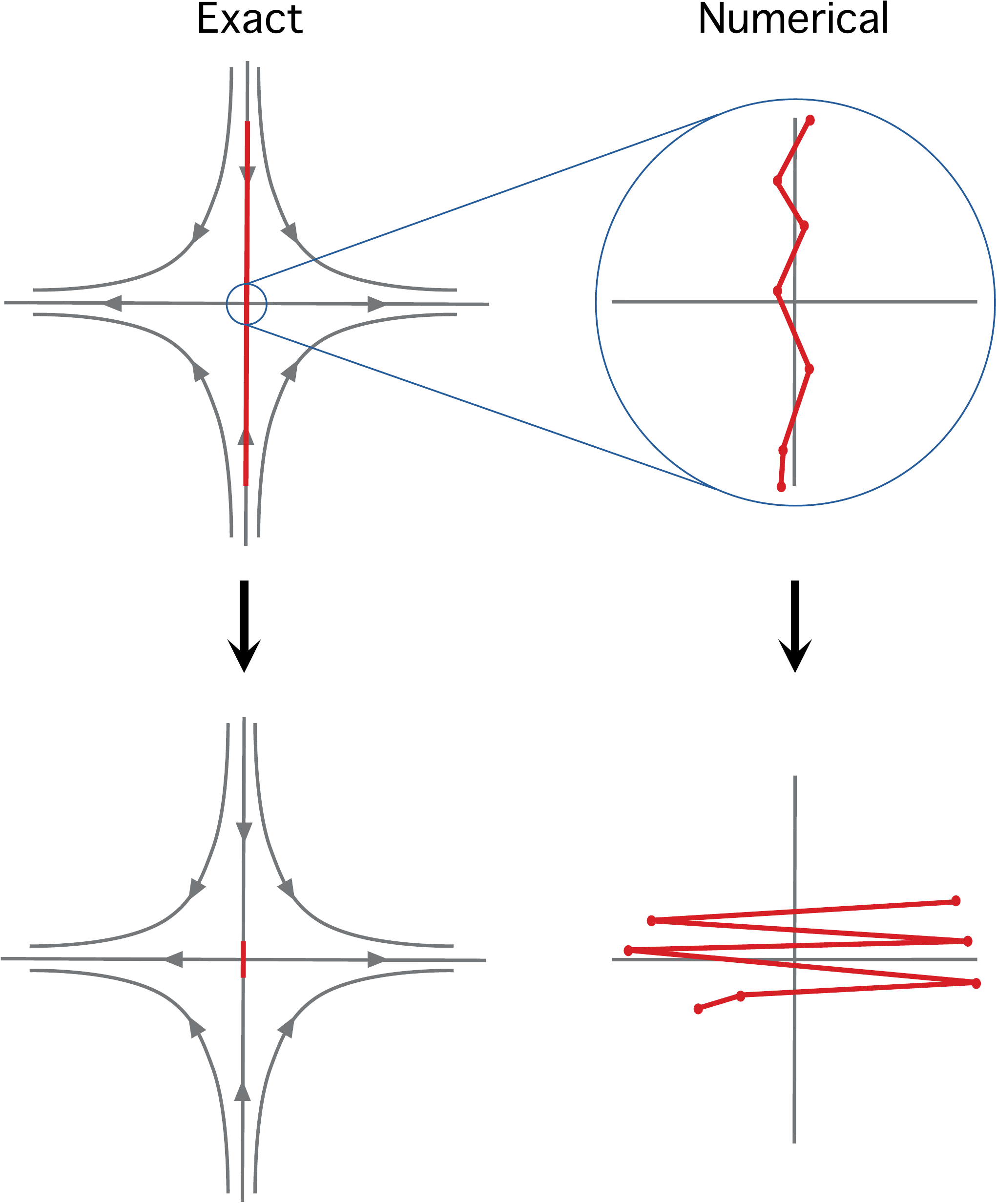}
\caption{Forward advection of a classic stable manifold (a repelling LCS over
finite times). Inaccuracies in determining the initial location of
the LCS lead to exponentially growing errors and accumulation along
the unstable manifold even if numerical errors were fully absent in
the advection.}
\label{fig:lambda-lemma}
\end{figure}

Most computational approaches to LCS seek their initial or final positions
as curves of initial conditions that lead to locally maximal trajectory
separation in forward or backward time. This repulsion-based approach
requires two numerical runs: one forward-time run that renders the
time-$t_{1}$ position of forward-repelling LCS, and one backward-time
run that reveals the time-$t_{2}$ position of forward-attracting
LCS. Determining the positions of these material surfaces at an intermediate
time $t$ accurately, however, comes at high computational cost: it
requires the accurate numerical advection of curves that are unstable
in the time direction of advection (see Fig.\ \ref{fig:lambda-lemma},
as well as the discussion in \cite{Farazmand2013}).

A recent computational advance is offered by \cite{Farazmand2013},
showing how both repelling and attracting LCS can be simultaneously
obtained either at $t_{1}$ or $t_{2}$ from a single numerical run.
This approach renders an attracting LCS at a time $t\in[t_{1},t_{2}]$
as the advected image of the initial LCS position at time $t_{1}$.
Similarly, the time-$t$ position of a repelling LCS can be obtained
by backward-advecting its position from time $t_{2}$ to $t$. Both
of these computations track attracting material surfaces, and hence
are numerically robust. However, they involve the advection of LCS
from two different initial times, and hence are necessarily preceded
by two separate numerical advections of a dense enough grid of initial
conditions. Altogether, therefore, the computational cost of constructing
both repelling and attracting LCS at arbitrary times $t\in[t_{1},t_{2}]$
has remained relatively high.

Here we propose a new computational strategy for two-dimensional incompressible
flows. Our strategy builds on results from \cite{Farazmand2013},
\cite{Karrasch2013d} and \cite{Schindler2012,Onu2013}, enabling
the reconstruction of all hyperbolic LCS for arbitrary times $t\in[t_{1},t_{2}]$
in a numerically robust fashion. This approach involves a single integration
of trajectories from a full numerical grid, followed by the advection
of select attracting material segments from local extrema of the singular
value field of the flow gradient. This procedure yields substantial
savings in computational time, as well as increased numerical accuracy
in LCS detection and tracking. We demonstrate these advantages on
a simple analytical flow example and on a direct numerical simulation
of two-dimensional turbulence.

This paper is organized as follows. In Section \ref{sec:Notation},
we fix our notation and recall relevant findings from \cite{Karrasch2013d}
on the singular value decomposition of the linearized flow map. In
Section \ref{sec:Attraction-LCS-Theory} we present our attraction-based
approach to hyperbolic LCS in the context of the recent geodesic theory
of LCS \cite{Haller2012,Farazmand2012-1,Farazmand2013a}. In Section \ref{sec:Examples},
we provide a proof of concept in the autonomous Duffing oscillator
and compare our approach to previous ones in a simulation of two-dimensional
turbulence, before concluding in Section \ref{sec:Conclusions}.

\section{Set-up\label{sec:Notation}}

Consider a smooth, two-dimensional vector field $v(x,t)$, defined
over a finite interval $I\coloneqq\left[t_{1},t_{2}\right]$ and over
spatial locations $x\in D\subset\mathbb{R}^{2}$. The trajectories generated
by $v(x,t)$ satisfy the ordinary differential equation
\begin{equation}
\dot{x}=v(x,t).\label{eq:ODE}
\end{equation}
The $t_{1}$-based flow map is denoted by $F_{t_{1}}^{t_{2}}\colon x_{1}\mapsto x_{2}$,
mapping initial values $x_{1}$ from time $t_{1}$ to their position
at time $t_{2}$ along the corresponding solution of \eqref{eq:ODE}.
We recall that the flow map is as smooth in $x_{1}$ as is $v$ in
$x$.

At any $x_{1}\in D,$ the \emph{deformation gradient }$DF_{t_{1}}^{t_{2}}(x_{1})$
is a matrix that admits a singular value decomposition (SVD)
of the form
\begin{align}
DF_{t_{1}}^{t_{2}} & =\Theta\Sigma\Xi^{\top}, & \Theta=\begin{pmatrix}\theta_{2} & \theta_{1}\end{pmatrix},\Xi=\begin{pmatrix}\xi_{2} & \xi_{1}\end{pmatrix} & \in O(2), & \Sigma & =\begin{pmatrix}\sigma_{2}^{\text{f}} & 0\\
0 & \sigma_{1}^{\text{f}}
\end{pmatrix},\label{eq:SVD}
\end{align}
with $\sigma_{2}^{\text{f}}\geq\sigma_{1}^{\text{f}}>0$ on the flow
domain $D$. The numbers $\sigma_{2}^{\text{f}},\sigma_{1}^{\text{f}}$
are the singular values\emph{ }of $DF_{t_{1}}^{t_{2}}$; the columns
of $\Xi$ (i.e., $\xi_{2}$ and $\xi_{1}$) are\emph{ }the right singular
vectors of $DF_{t_{1}}^{t_{2}}$; the columns of $\Theta$ (i.e.,
$\theta_{2}$ and $\theta_{1}$) are the left singular vectors\emph{
}of $DF_{t_{1}}^{t_{2}}$. From \eqref{eq:SVD}, we see that
\begin{align*}
DF_{t_{1}}^{t_{2}}\left(x_{1}\right)\xi_{i}\left(x_{1}\right) & =\sigma_{i}^{\text{f}}\left(x_{1}\right)\theta_{i}\left(x_{2}\right), & x_{2} & =F_{t_{1}}^{t_{2}}(x_{1}), & i & \in\left\{ 1,2\right\} .
\end{align*}
We recall that the singular values $\sigma_{2}^{\text{f}}(x_{1})$
and $\sigma_{1}^{\text{f}}(x_{1})$ measure infinitesimal stretching
and compression along the trajectory starting from $x_{1}$. Furthermore,
the unit vectors $\xi_{2}\left(x_{1}\right)$ and $\xi_{1}\left(x_{1}\right)$
are the tangent vectors pointing to the directions of strongest stretching
and compression under the linearized flow $DF_{t_{1}}^{t_{2}}(x_{1})$.

If the velocity field is incompressible, i.e., $\nabla_{x}\cdot v(x,t)\equiv0$,
then $\det\left(DF\right)=\sigma_{1}^{\text{f}}\sigma_{2}^{\text{f}}=1$,
and consequently
\begin{align}
\sigma_{2}^{\text{f}} & =1/\sigma_{1}^{\text{f}}.\label{eq:svd_incompressible}
\end{align}
As a result, local maxima of $\sigma_{2}^{\text{f}}$ (locally strongest-stretching
points) coincide with local minima of $\sigma_{1}^{\text{f}}$ (locally
strongest-compressing points). At any point $x_{1}\in D$, the average
exponential rate of largest stretching over the time interval $\left[t_{1},t_{2}\right]$
of length $T=t_{2}-t_{1}$ is defined as
\[
\Lambda^{\text{f}}\left(x_{1}\right)\coloneqq\frac{1}{T}\log\sigma_{2}^{\text{f}}\left(x_{1}\right),
\]
which is referred to as the \emph{(forward) finite-time Lyapunov exponent}
\emph{(FTLE)}. In the incompressible case, Eq.\ \eqref{eq:svd_incompressible}
shows that the FTLE can equally well be considered as a measure of
the strongest local compression at $x_{1}$.

For the backward flow from $t_{2}$ to $t_{1}$, the backward deformation
gradient is given by
\begin{equation}
DF_{t_{2}}^{t_{1}}(x_{2})=\left[DF_{t_{1}}^{t_{2}}(x_{1})\right]^{-1}=\Xi\Sigma^{-1}\Theta^{\top},\label{eq:svd_inverse}
\end{equation}
with $\Sigma^{-1}=\begin{pmatrix}\sigma_{1}^{\text{b}} & 0\\
0 & \sigma_{2}^{\text{b}}
\end{pmatrix}=\begin{pmatrix}1/\sigma_{2}^{\text{f}} & 0\\
0 & 1/\sigma_{1}^{\text{f}}
\end{pmatrix}$. The singular values of $DF_{t_{2}}^{t_{1}}$ are therefore given
by
\begin{equation}
\sigma_{2}^{\text{b}}(x_{2})=1/\sigma_{1}^{\text{f}}(x_{1}),\qquad\sigma_{1}^{\text{b}}(x_{2})=1/\sigma_{2}^{\text{f}}(x_{1}),\qquad x_{2}=F_{t_{1}}^{t_{2}}(x_{1}),\label{eq:svd_identities}
\end{equation}
and the backward right singular vectors are given by $\theta_{1}$
and $\theta_{2}$, the strongest- and weakest-stretching directions
at $x_{2}$ in backward time.%
\footnote{The superscripts $\text{f}$ and $\text{b}$ refer to forward and
backward time, respectively. %
}

Eq.\ \eqref{eq:svd_incompressible} shows that the maximal (minimal)
singular value of the linearized flow map is equal to the maximal
(minimal) singular value of the linearized inverse flow map. Thus,
local maxima of $\sigma_{2}^{\text{f}}$ are mapped bijectively to
local maxima of $\sigma_{2}^{\text{b}}$ by the flow map. We summarize
the equivalences of local extrema of the forward and backward singular
value fields as follows:
\begin{equation}
\begin{array}{cccc}
 & x_{1}\text{ at }t_{1} & \xrightarrow{F_{t_{1}}^{t_{2}}} & x_{2}=F_{t_{1}}^{t_{2}}(x_{1})\text{ at }t_{2}\\
\\
 & \sigma_{2}^{\text{f}}\text{--maximum} & \Longleftrightarrow & \sigma_{1}^{\text{b}}\text{--minimum}\\
\text{if incompressible} & \Updownarrow &  & \Updownarrow\\
 & \sigma_{1}^{\text{f}}\text{--minimum} & \Longleftrightarrow & \sigma_{2}^{\text{b}}\text{--maximum}.
\end{array}\label{eq:equivalence}
\end{equation}
In terms of the backward FTLE field, we recover \cite[Prop. 2]{Haller2011-1}
for the incompressible case:
\begin{equation*}
\Lambda^{\text{b}}\left(x_{2}\right)=\frac{1}{T}\log\sigma_{2}^{\text{b}}\left(x_{2}\right)=\frac{1}{T}\log\left(\sigma_{1}^{\text{f}}\left(x_{1}\right)\right)^{-1}=\frac{1}{T}\log\sigma_{2}^{\text{f}}\left(x_{1}\right)=\Lambda^{\text{f}}\left(x_{1}\right).
\end{equation*}

In summary, as argued in \cite{Karrasch2013d}, the SVD of $DF_{t_{1}}^{t_{2}}$
yields complete forward and backward stretch information from a uni-directional
flow computation.

\section{Forward and backward Geodesic Theory of Hyperbolic LCS \label{sec:Attraction-LCS-Theory}}

The following definition recalls the hyperbolic LCS candidates obtained
from two-di\-men\-sion\-al geodesic LCS theory.

\begin{defn*}[Shrink and stretch lines, \cite{Farazmand2012-1,Farazmand2013}]
We call a smooth curve $\gamma$ a \emph{forward (}or \emph{backward)
shrink line,} if it is pointwise tangent to the $\xi_{1}$ (or $\theta_{2}$)
field. Similarly, we call $\gamma$ a \emph{forward (}or \emph{backward)
stretch line}, if it is pointwise tangent to the $\xi_{2}$ (or $\theta_{1}$)
field.
\end{defn*}

Shrink and stretch lines are solutions of a variational principle
put forward in \cite{Farazmand2013a} for LCS. This principle stipulates
as a necessary condition that the time $t_{1}$ positions of hyperbolic
LCS must be stationary curves of the averaged Lagrangian shear \cite{Farazmand2013a}.
This variational principle leads to the result that time $t_{1}$
positions of hyperbolic LCS are necessarily null-geodesics of an appropriate
Lorentzian metric associated with the deformation field \cite{Farazmand2013a}.
This prompts us to refer to the underlying approach as geodesic LCS
theory.

Away from points where $\sigma_{2}^{\text{f}}=\sigma_{1}^{\text{f}}$
at $t=t_{1}$ and $\sigma_{2}^{\text{b}}=\sigma_{1}^{\text{b}}$ at
$t=t_{2}$, both the initial and the final flow configuration is foliated
continuously by mutually orthogonal forward and backward shrink and
stretch lines. As discussed in \cite{Farazmand2013,Karrasch2013d},
the following equivalence relations hold:
\begin{equation}
\begin{array}{ccc}
\text{at }t_{1} & \xrightarrow{F_{t_{1}}^{t_{2}}} & \text{at }t_{2}\\
\text{}\\
\text{forward shrink line} & \Longleftrightarrow & \text{backward stretch line}\\
\perp &  & \perp\\
\text{forward stretch line} & \Longleftrightarrow & \text{backward shrink line}.
\end{array}\label{eq:equivalence_2}
\end{equation}

The forward shrink and stretch lines provide candidate curves for
the positions of repelling and attracting LCS at time $t_{1}$. To
find the positions of actual hyperbolic LCS as centerpieces of observed
tracer deformation, we seek members of these two line families that
evolve into locally most attracting or repelling material lines over
the time interval $[t_{1},t_{2}].$

To this end, we follow \cite{Schindler2012,Onu2013} to require a
sufficient condition that hyperbolic LCS must satisfy. Specifically,
the time $t_{1}$ positions of \emph{forward repelling LCS} are shrink
lines that intersect local maxima of $\sigma_{2}^{\text{f}}$\emph{;
}the time $t_{1}$ positions of\emph{ forward attracting LCS} are
stretch lines that intersect local minima of $\sigma_{1}^{\text{f}}$.
The time $t_{2}$ positions of backward-repelling and backward-attracting
LCS are defined analogously using the backward singular-value fields
$\sigma_{2}^{\text{b}}$ and $\sigma_{1}^{\text{b}}$. By the equivalences
detailed above, forward-attracting LCS, as evolving material lines,
coincide with backward repelling LCS. Similarly, forward-repelling
LCS, as evolving material lines, coincide with backward-attracting
LCS.

If the vector field $v(x,t)$ is incompressible, then the relation
\eqref{eq:svd_incompressible} forces local maxima of $\sigma_{2}^{\text{f}}$\emph{
}to coincide with local minima of $\sigma_{1}^{\text{f}}$. As a consequence,
both forward-repelling and forward-attracting LCS intersect the maxima
of $\sigma_{2}^{\text{f}}$ at time $t_{1}$. This fact will simplify
our upcoming computational algorithm considerably for incompressible
flows.

As noted earlier, reconstructing a full forward-attracting LCS as
a material line involves advecting its time $t_{1}$ position under
the flow map. This is a self-stabilizing numerical process, as it
tracks an attracting surface. In contrast, reconstructing a forward-repelling
LCS from its time $t_{1}$ position by flow advection is an unstable
numerical process. Indeed, the smallest initial errors in identifying
the LCS position are quickly amplified, as shown in Fig.\ \ref{fig:lambda-lemma}.

Relations \eqref{eq:equivalence} and \eqref{eq:equivalence_2}, however,
allow us to compute the forward-repelling LCS equivalently as backward-attracting
LCS. Specifically, forward-repelling LCS positions at a time $t\in[t_{1},t_{2}]$
can be equivalently obtained from advection under the backward flow
map $F_{t_{2}}^{t}.$ The curves to be advected under $F_{t_{2}}^{t}$
are just the backward stretch lines running through local minima of
$\sigma_{1}^{\text{b}}$. By \eqref{eq:svd_identities}, however,
local minima of $\sigma_{1}^{\text{b}}$ are just the images of local
maxima of $\sigma_{2}^{\text{f}}$ under the flow map $F_{t_{1}}^{t_{2}}$.

\begin{figure}
\centering
\includegraphics[width=0.4\textwidth]{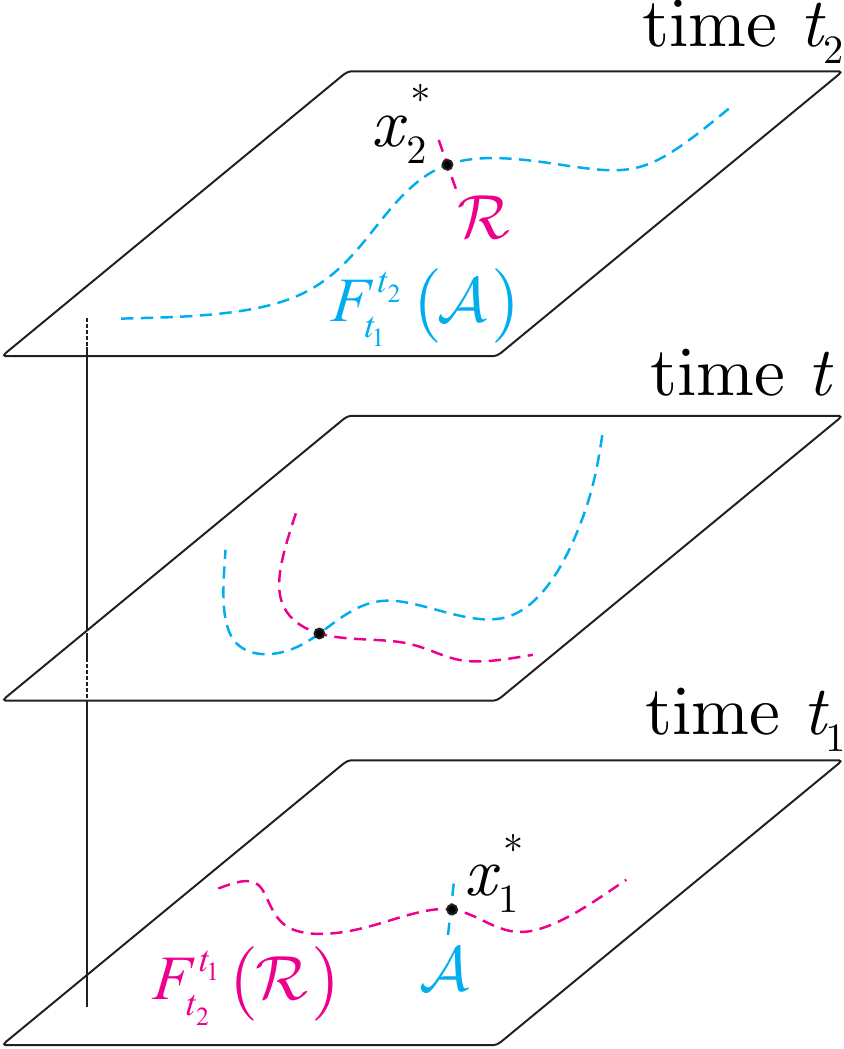}
\caption{Illustration of the attraction-based LCS extraction for an incompressible
flow at an arbitrary time $t\in[t_{1},t_{2}].$ Here $\mathcal{A}$
denotes a short vector parallel to $\xi_{2}$ at a local maximum $x_{1}^{*}$
of the $\sigma_{2}^{\text{f}}(x)$ field . Similarly, $\mathcal{R}$
denotes a short vector parallel to $\theta_{1}$ at the point $x_{2}^{*}=F_{t_{1}}^{t_{2}}(x_{1}^{*})$.
Recall that both forward-repelling and forward-attracting LCS intersect
the maxima of $\sigma_{2}^{\text{f}}$ at time $t_{1}$ in case of
incompressibility.}
\label{fig:new_approach}
\end{figure}

The computation of stretch lines still involves the integration of
direction fields, for which orientation issues have to be resolved
(see \cite{Tchon2006,Farazmand2012-1}). A new feature we introduce
here is to advect short line segments (tangents) as opposed to whole
stretch lines running through the appropriate extrema of the singular
value fields. This idea exploits the tangentially stretching and normally
attracting nature of stretch lines, saves on computational cost, and
produces highly accurate results, as we demonstrate later. We summarize
our attraction-based LCS construction in Fig.\ \ref{fig:new_approach}
for the case of incompressible flows. For compressible flows, forward-attracting
LCS at time $t_{1}$ are still constructed from local minima of $\sigma_{1}^{\text{f}}$,
but backward-attracting LCS at $t_{2}$ are constructed from advected
local maxima of $\sigma_{2}^{\text{f}}$, which generally differ from
advected local minima of $\sigma_{1}^{\text{f}}$.

\subsection*{Numerical implementation}

Here we summarize the computational steps resulting from our previous
considerations, assuming a forward-time advection of the chosen numerical
grid.
\begin{enumerate}
\item \textbf{Compute flow map and its linearization:} We solve the ODE
\eqref{eq:ODE} from a sufficiently dense grid of initial conditions
to obtain a discrete approximation to the flow map $F_{t_{1}}^{t_{2}}$.
We also obtain a numerical approximation to the linearized flow map
$DF_{t_{1}}^{t_{2}}$ at the grid points by one of four methods: (i)
solving the equation of variations associated with \eqref{eq:ODE},
(ii) finite-differencing $F_{t_{1}}^{t_{2}}$ along the grid, (iii)
finite-differencing on a smaller auxiliary grid \cite{Farazmand2012-1},
(iv) via convolution with Gaussian kernels \cite{Peikert2008}.
\item \textbf{Compute singular values: }We compute the singular-value decomposition
of the deformation gradient tensor field $DF_{t_{1}}^{t_{2}}$. This
yields the singular values $\sigma_{i}^{\text{f }}$ as well as the
right- and left-singular vector fields $\xi_{i}$ and $\theta_{i}$,
respectively. The singular values $\sigma_{i}^{\text{b }}$ are obtained
directly from the relation \eqref{eq:svd_identities}.
\item \textbf{Select seeding points for LCS:} We need to identify points
of strongest attraction, i.e. local minima of $\sigma_{1}^{\text{f}}$
at the initial time and local minima of $\sigma_{1}^{\text{b}}$ at
the final time. While the first are identified directly, the latter
are advected images of local maxima of $\sigma_{2}^{\text{f}}$ under
the flow map $F_{t_{1}}^{t_{2}}$. In the incompressible case, the
points of strongest attraction coincide with local maxima of $\sigma_{2}^{\text{f}}$ and
their advected images under $F_{t_{1}}^{t_{2}}$, respectively. As
in \cite{Onu2013}, we start by sorting all local maxima in ascending
order by the values of $\sigma_{1}^{\text{f}}$ or descending order
by the values of $\sigma_{2}^{\text{f}}$. We then pick the first
point $p_{1}$ in the ordered list and discard all local extrema in
a small neighborhood of $p_{1}$. From the remaining points on the
list, we pick the first point $p_{2}$ and discard extrema in a small
neighborhood of $p_{2}$, and so on. This procedure filters out local
extrema in noisy singular value fields.
\item \textbf{Compute hyperbolic LCS}: For any time $t\in[t_{1},t_{2}]$
of interest, we use the flow map $F_{t_{1}}^{t}$ to advect short
line segments tangent to $\xi_{2}(p_{i})$ at the points $p_{i}$
identified in the previous step. The resulting set of curves form
the time $t$ positions of attracting LCS. In the incompressible case,
we use the flow map $F_{t_{2}}^{t}$ to advect short line segments
tangent to $\theta_{1}(F_{t_{2}}^{t}(p_{i}))$ at the points $F_{t_{1}}^{t_{2}}(p_{i})$.
Recall that the characteristic stretching directions for the backward flow are
obtained from the forward time computation in step (2) due to Eq.\ \eqref{eq:svd_inverse}.
The resulting set of curves form the time $t$ positions of repelling
LCS. For the advection of line segments, the use of an adaptive integration
scheme may be necessary. This is to fill emerging large gaps between
adjacent points due to stretching, and to mitigate the possibly high
curvature in the tracked material curve (see, e.g., \cite{Mancho2003}).
\end{enumerate}

\section{Examples\label{sec:Examples}}

\subsection{Duffing oscillator}

We first consider a rescaled version of the unforced, undamped Duffing
oscillator with Hamiltonian
\[
H(x,y)=\frac{1}{4}x^{4}-2x^{2}+\frac{1}{2}y^{2}.
\]
This example has already been used to illustrate shrink and stretch
line context by \cite{Farazmand2013} locally around the origin, showing
the convergence of forward- and backward maximal stretch directions
to the unstable and stable subspaces, respectively. In our present
computations, we use the times $t_{1}=0$ and $t_{2}=T=2.5$.

\begin{figure}
\centering
\includegraphics[height=0.18\textheight]{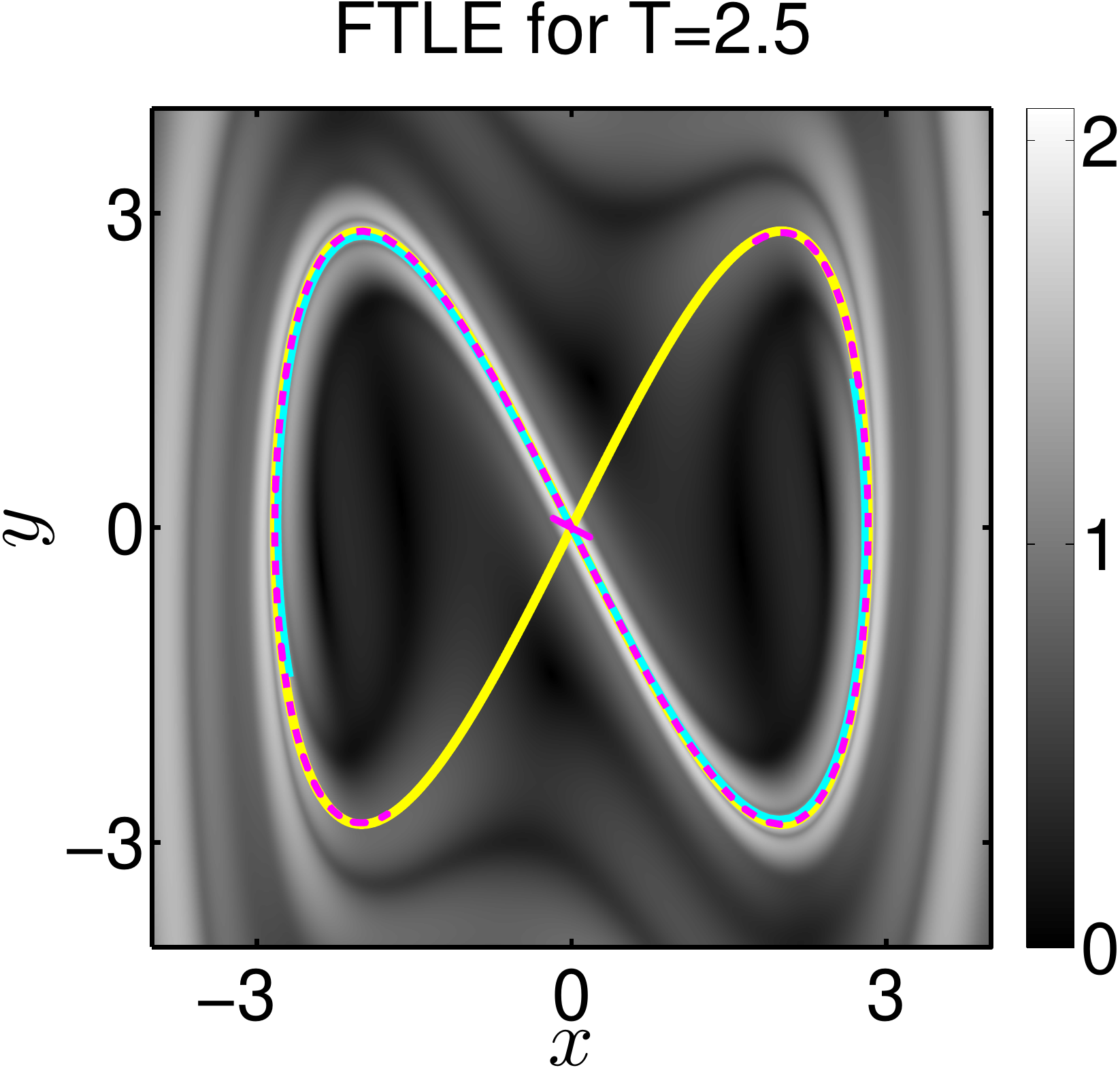}~%
\includegraphics[height=0.18\textheight]{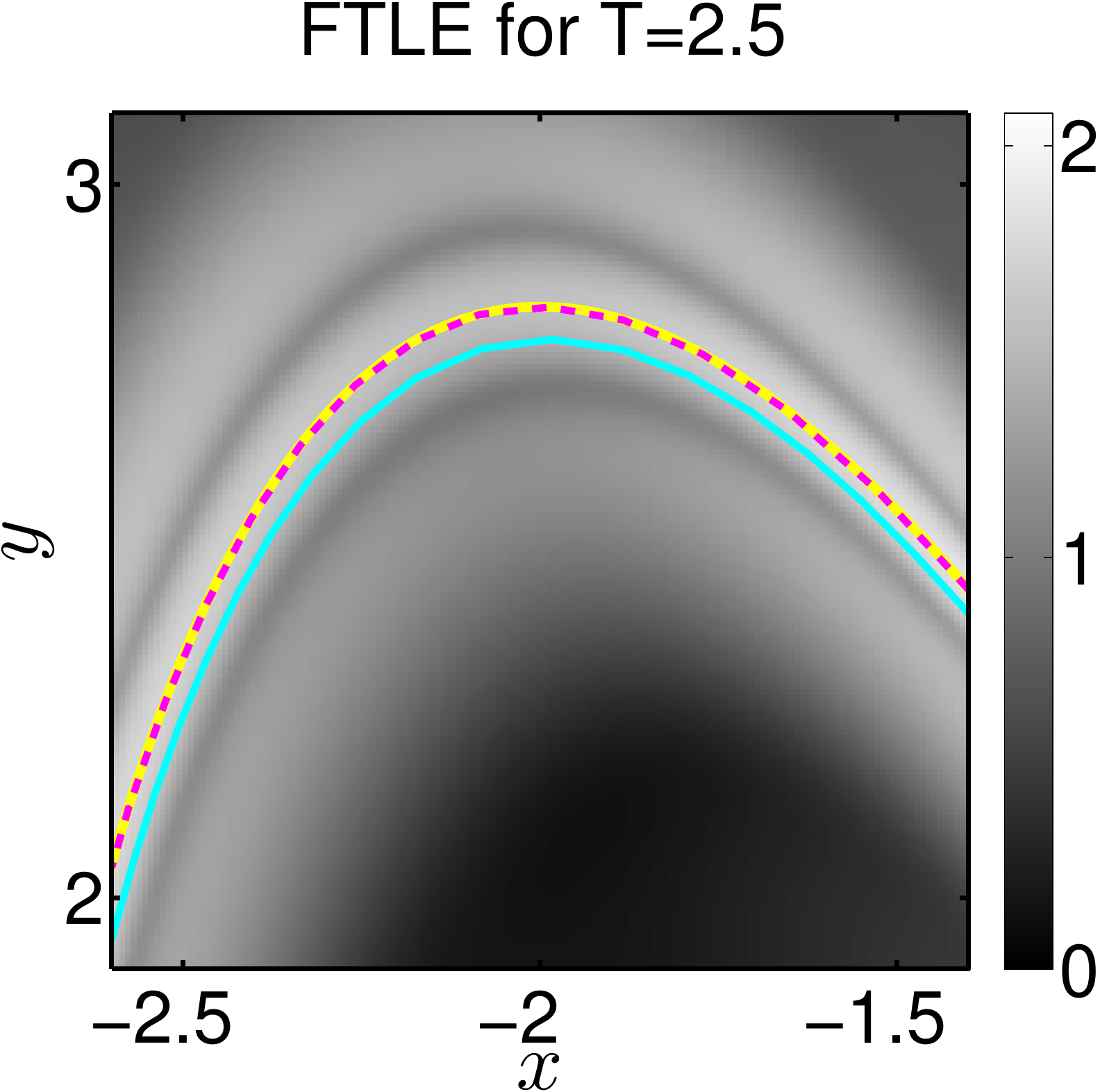}~%
\includegraphics[height=0.18\textheight]{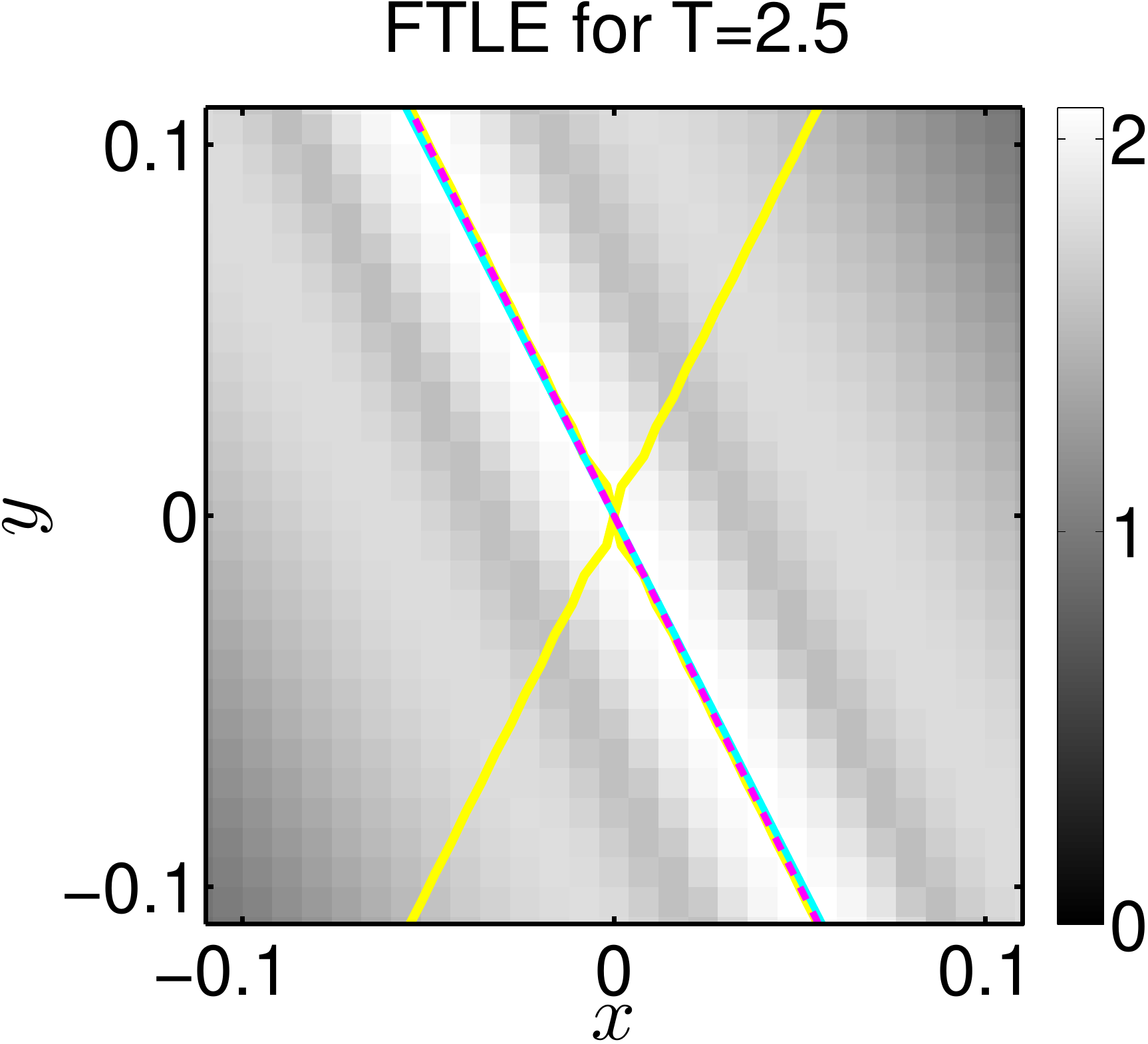}
\caption{In the background, the FTLE field for integration time $T=2.5$ is
shown. On top, we show the zero energy level $H=0$ (yellow), the
shrink line (cyan), and a straight line aligned with $\theta_{1}(0)$
(short line segment not aligned with the homoclinic, magenta) at $t_{2}=2.5$
together with its image at $t_{1}=0$ (magenta). The left figure shows
that all structures highlight the homoclinic loop with reasonable
accuracy. The magnification in the middle shows, however, a significant
deviation of the shrink line from the stable manifold. At the same
time, the backward-advected straight line segment approximates the
stable manifold perfectly. The right plot shows both the shrink line
and the image of the backward stretch line segment to perform well
near the origin.}
\label{fig:duffing_2}
\end{figure}

In Fig.\ \ref{fig:duffing_2}, we compare the results from the earlier
numerical LCS detection scheme used in \cite{Haller2012} to our approach
described in Section \ref{sec:Attraction-LCS-Theory}. While the left
plot shows all structures to highlight the homoclinic loop, the middle
plot shows that the shrink line deviates from the loop visibly at
the first turn. In contrast, the backward-advected line segment stays
close to the loop. The right plot shows that at the origin, both the
shrink line and the advected stretch line indicate consistently the
direction of strongest attraction.

\begin{figure}
\centering
\includegraphics[height=0.2\textheight]{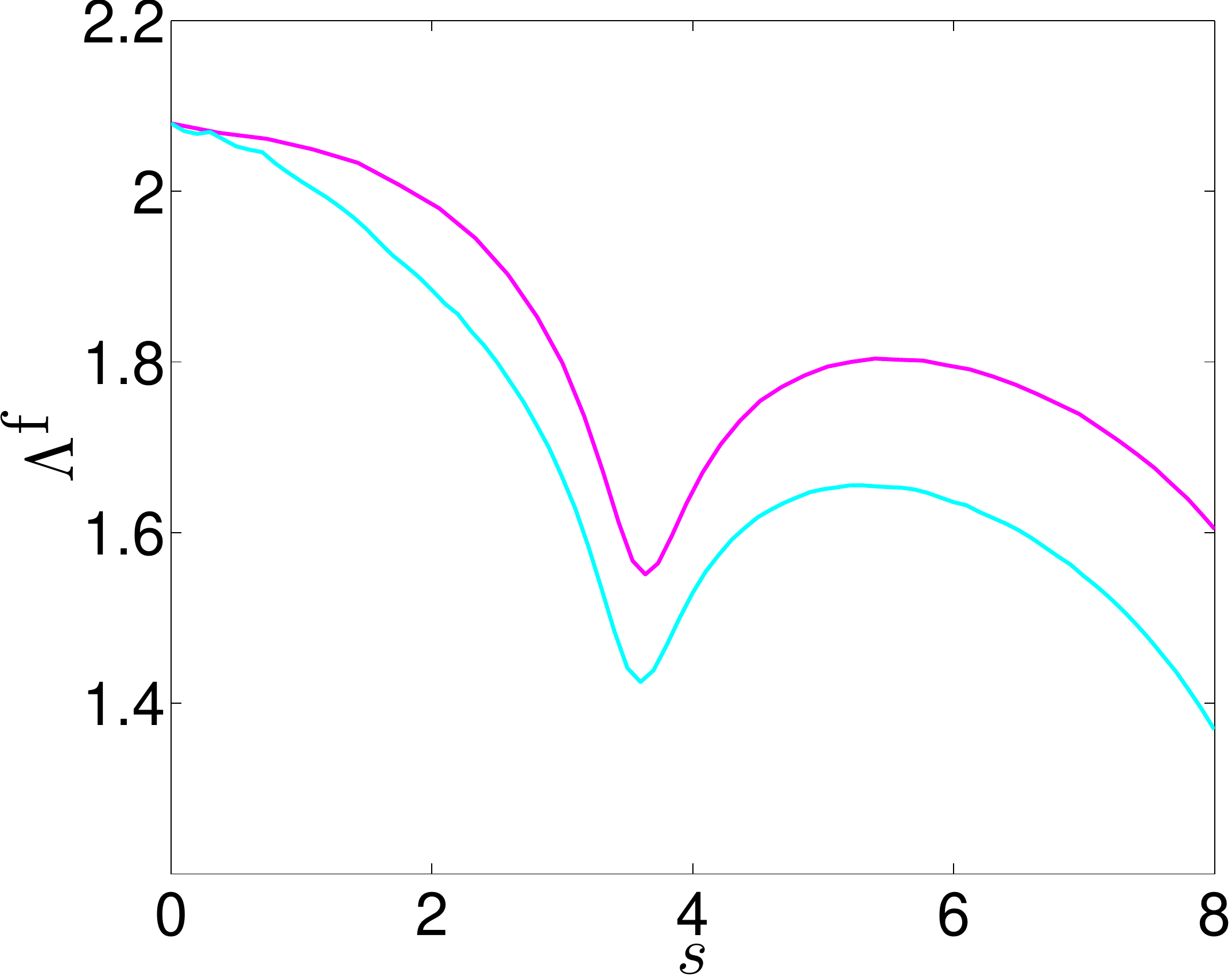}~%
\includegraphics[height=0.22\textheight]{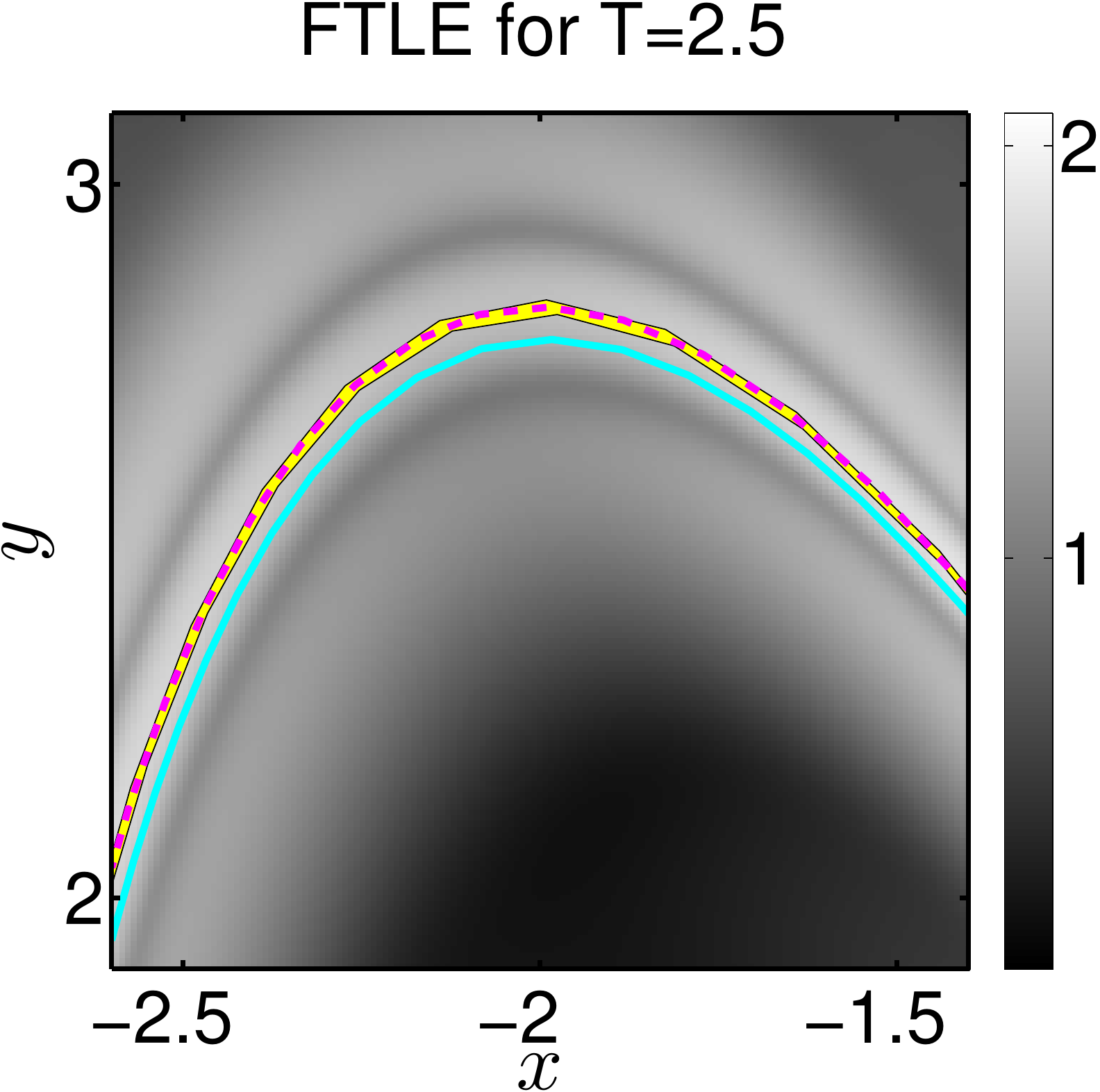}
\caption{Left: Comparison of FTLE along the backward-advected backward stretch
line (magenta) and the forward shrink line (cyan). Note that the advected
stretch line has a uniformly higher repulsion rate and is therefore
a better approximation to the repelling LCS. Right: Backward-advected
particle blob of initial diameter $1.0$ (yellow), backward-advected
stretch line (dashed magenta) and forward shrink line (cyan), showing
that the advected stretch line is a better approximation to the backward
attracting (i.e.\ repelling) LCS. (The numerical advection is performed
by the \textsc{Matlab} routine \texttt{ode45} with absolute and relative
error tolerance of $10^{-8}$.)}
\label{fig:hyp_test}
\end{figure}

Fig.\ \ref{fig:hyp_test} gives further quantitative evidence that
the backward-advected backward stretch line gives a better approximation
to the actual repelling LCS position at time $t_{1}$ than the direct
computation of this LCS position from forward shrink lines.

Even in this simple example, therefore, the actual evolution of a
shrink line and a backward-advected backward stretch line are noticeable
different, although they should theoretically be identical. The root cause
is numerical errors in the singular vector computation,
as well as the limited ability of the discrete numerical grid to approximate
a repelling LCS (local stable manifold) as a continuous curve. The
error is initially invisible, but starts to accumulate rapidly during
integration of the $\xi_{2}$ ($\theta_{1}$) field and advection.

\subsection{Two-dimensional turbulence}

As a second example, we consider the two-dimensional Navier--Stokes
equations
\begin{align*}
\partial_{t}v+v\cdot\nabla v & =-\nabla p+\nu\Delta v+f,\\
\nabla\cdot v & =0,\\
v(\cdot,0) & =v_{0},
\end{align*}
where the unsteady velocity field $v(x,t)$ is defined on the two-dimensional
domain $U=[0,2\pi]\times[0,2\pi]$ with doubly periodic boundary conditions.
As in \cite{Farazmand2013,Farazmand2014}, we use a standard pseudo-spectral
method with 512 modes in each direction, and $2/3$ de-aliasing to
solve the above Navier--Stokes equations with viscosity $\nu=10^{-5}$
on the time interval $\left[0,100\right]$. The flow integration is
then carried out over the interval $t\in[50,100]$, in which the turbulent
flow has fully developed, by a fourth-order Runge--Kutta method with
variable step-size. The initial condition $v_{0}$ is the instantaneous
velocity field of a decaying turbulent flow. The external force $f$
is random in phase and band-limited, acting on the wave-numbers $3.5<k<4.5$.

\begin{figure}
\centering
\includegraphics[width=0.3\textwidth]{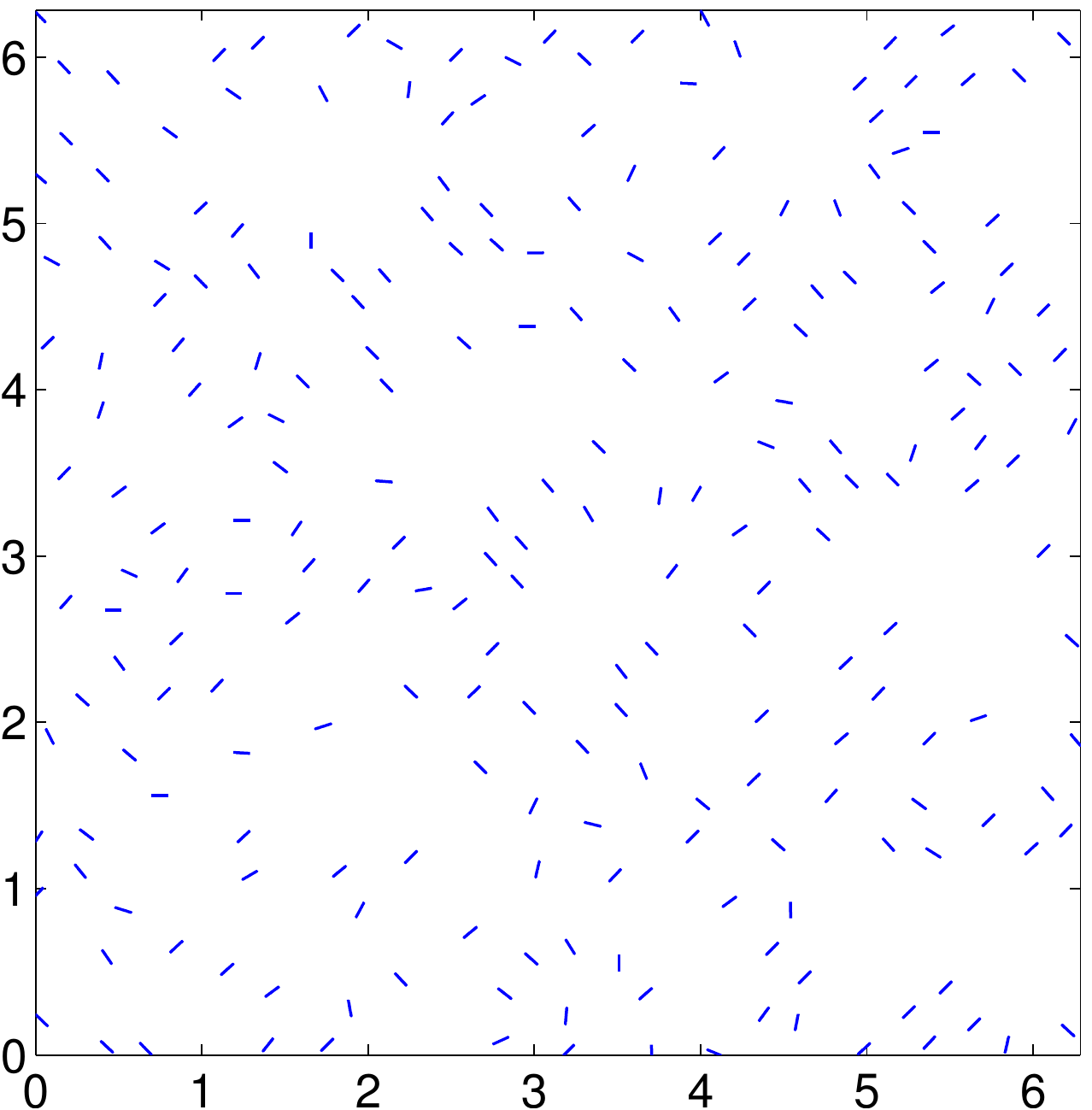}~~%
\includegraphics[width=0.3\textwidth]{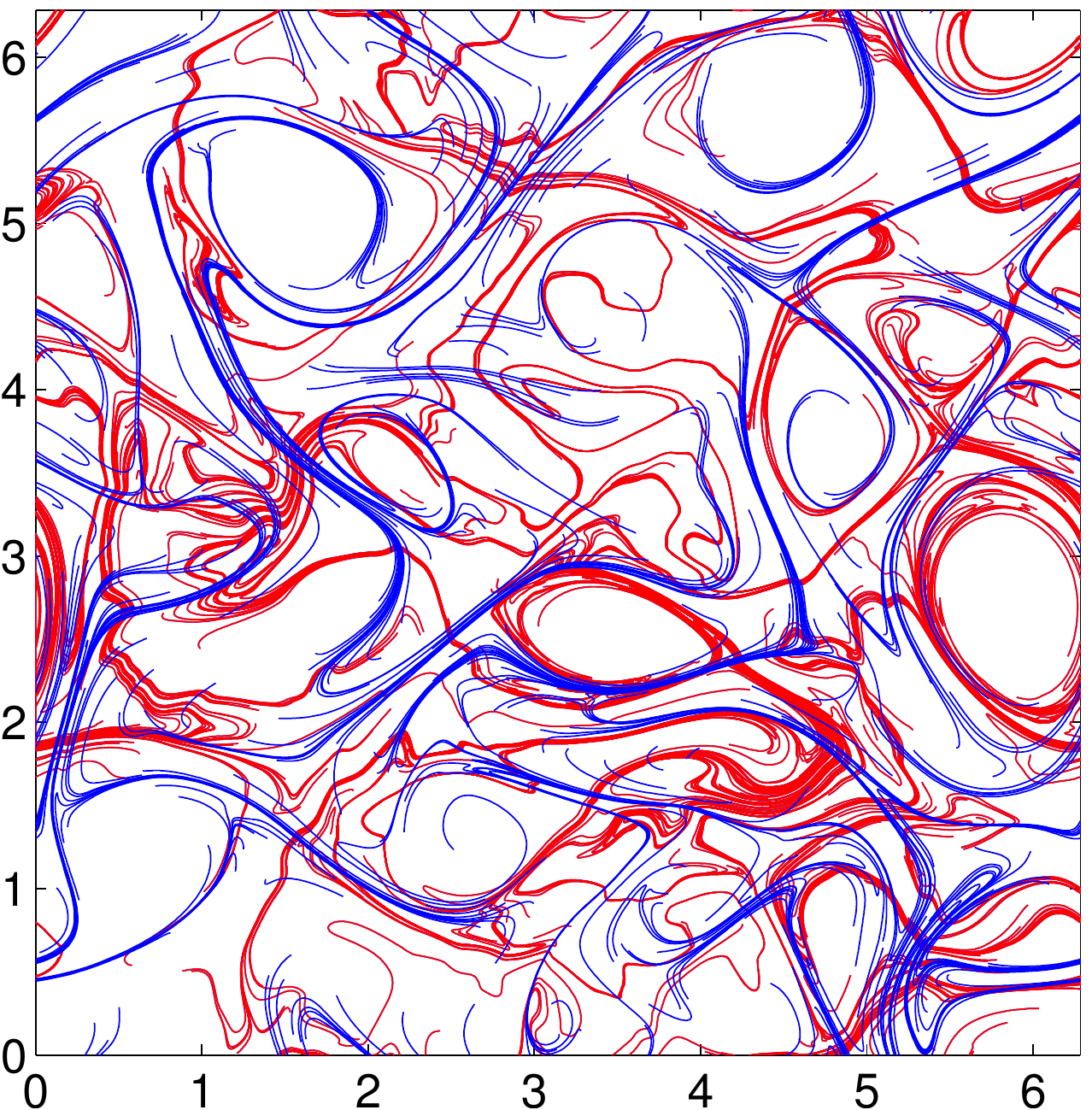}~~%
\includegraphics[width=0.3\textwidth]{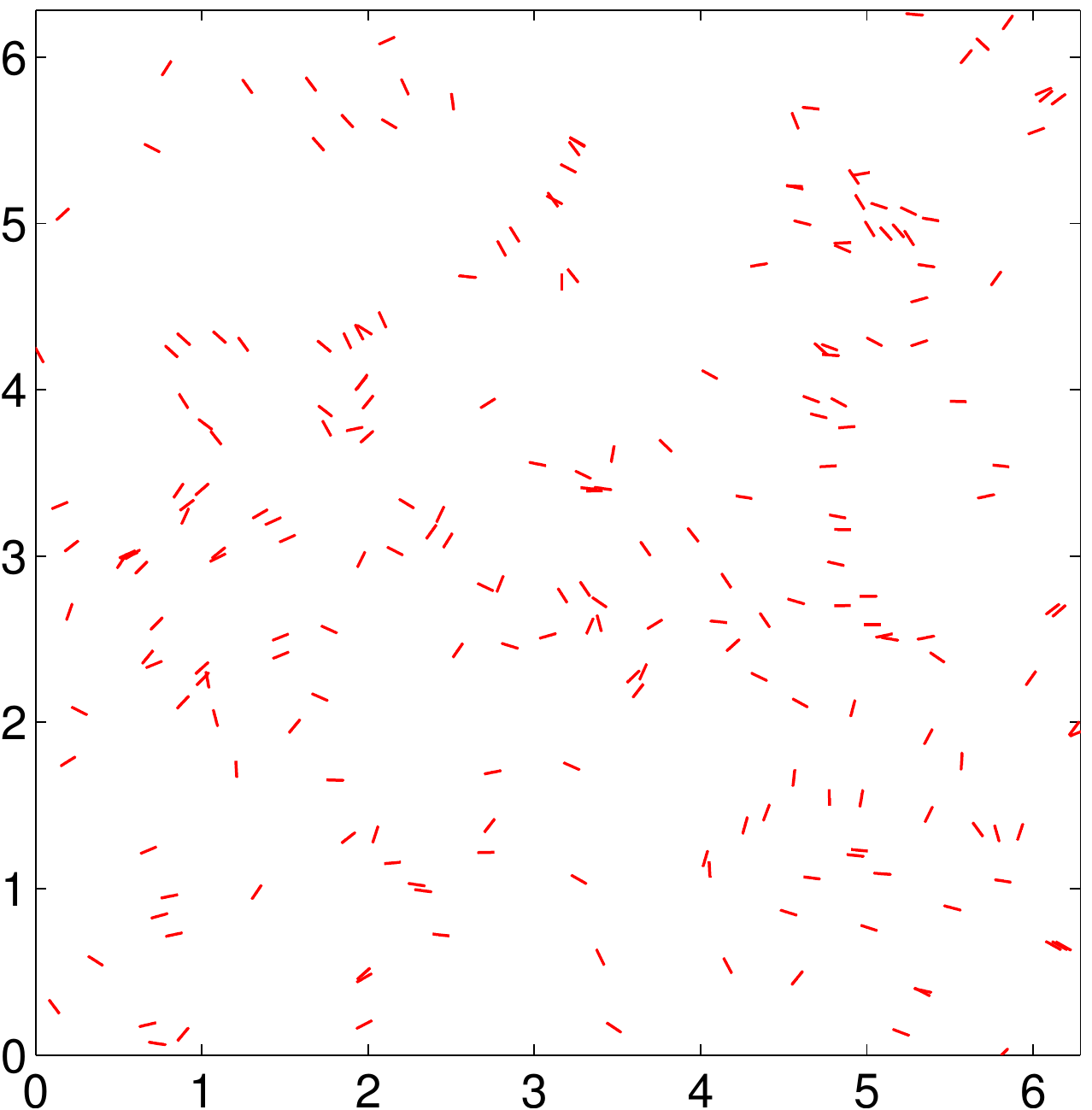}
\caption{Attracting (blue) and repelling (red) LCS in a simulation of two-dimensional
turbulence over the time interval $\left[50,100\right]$. Left: Initial
line segments at $t_{1}=50$ for the attracting LCS. Middle: Hyperbolic
LCS positions at $t=75$. Right: Initial line segments at $t_{2}=100$
for the repelling LCS.}
\label{fig:turbulence}
\end{figure}

\begin{figure}
\centering
\includegraphics[width=0.95\textwidth]{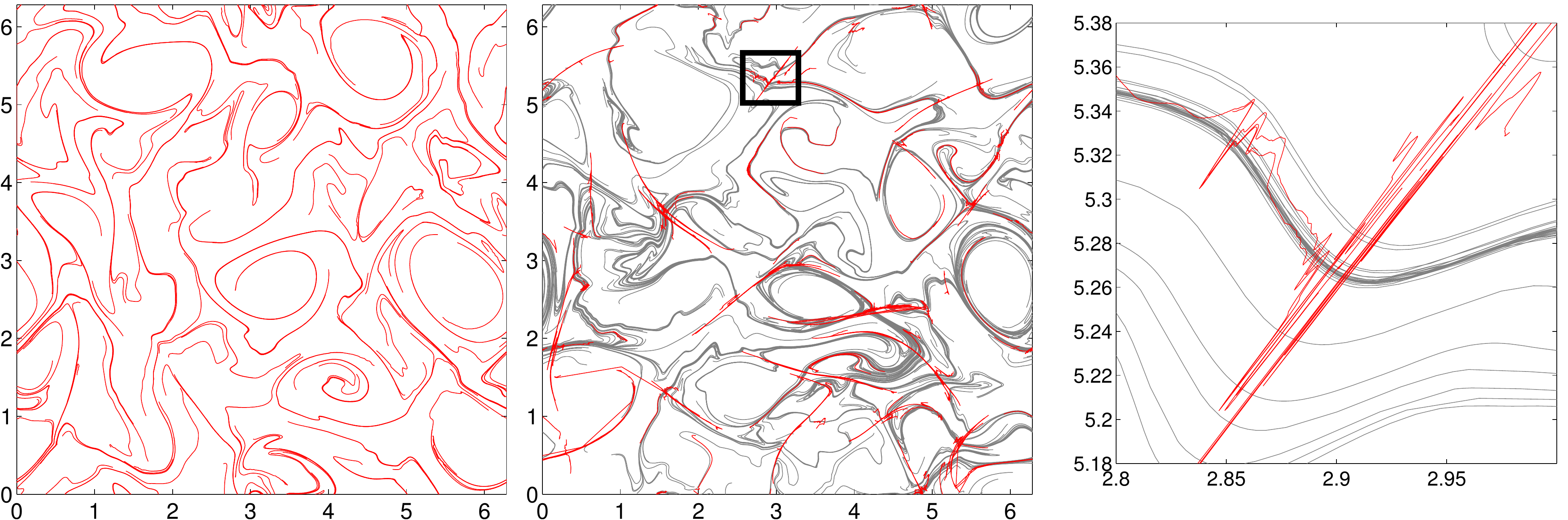}
\caption{Left: shrink lines computed directly at $t_{1}=50$ as curves tangent
to the $\xi_{1}(x)$ line field that intersect local maxima of $\sigma_{2}^{\text{f}}$.
Middle: the same shrink lines (in red) advected to $t=75$ to highlight
repelling LCS positions at that time. The gray curves are backward-advected
stretch lines from $t_{2}=100$ that run through the time $t_{2}$
positions of trajectories starting from local maxima of $\sigma_{2}^{\text{f}}$
at time $t_{1}$. Right: a close-up view of the middle panel, clearly
showing dramatic local inaccuracies from the forward calculation,
resembling the effect shown in Fig.\ \ref{fig:lambda-lemma}.}
\label{fig:turbulence_2}
\end{figure}

In Fig.\ \ref{fig:turbulence}(middle), we plot repelling (red) and
attracting (blue) LCS at the middle time instance $t=75$. As described
in Section \ref{sec:Attraction-LCS-Theory}, these LCS were launched
as straight line segments of length $0.1$ from local $\sigma_{2}^{\text{f}}$--maxima
and their flow images, which are $\sigma_{2}^{\text{b}}$--maxima,
see Fig.\ \ref{fig:turbulence}(left) and (right). The filtering
radius for local $\sigma_{2}^{\text{f}}$--maxima was set to $0.2$,
yielding a reduction from $11,000$ to $229$ seeding points.

We plot forward shrink lines at the initial time $t_{1}=50$ in Fig.\ \ref{fig:turbulence_2}(left),
and compare their forward-advected images (red) at the intermediate
time $t=75$ with the backward advected stretch lines (gray), seeded
at the corresponding points (see the middle panel of Fig.\ \ref{fig:turbulence_2}).
Analytically, these curves should coincide. In some locations, they
indeed agree well, but in other locations, the discrepancy is dramatic
(see the close-up view in the right panel of Fig.\ \ref{fig:turbulence_2}).
This is the consequence of the effect illustrated in Fig.\ \ref{fig:lambda-lemma},
showing the clear advantage of our method over the forward-time tracking
of a repelling LCS.

\section{Conclusions\label{sec:Conclusions}}

We have proposed a paradigm shift in the detection of hyperbolic Lagrangian
Coherent Structures (LCS). Instead of detecting initial positions
of LCS as curves of maximal forward repulsion, we seek them as backward-advected
locations of maximal backward attraction. While these two approaches
are theoretically equivalent, the latter approach (developed here)
eliminates an inherent numerical instability of the former approach
(used in prior work). We have demonstrated that our attraction-based
approach leads to substantial improvements in accuracy and computational
cost.

We have discussed our approach in the framework of the geodesic theory
\cite{Haller2012,Farazmand2013,Farazmand2013a}, because this theory
allows for the explicit computation of hyperbolic LCS as parametrized
curves. The proposed focus on attraction, however, automatically extends
to potential future refinements in LCS computations.

The advection of identified hyperbolic LCS in the stable time direction
is a simple idea, but relies heavily on the notion of a forward-time
attracting LCS, which has been proposed only recently \cite{Farazmand2013}.
We have combined this notion with the SVD of the deformation gradient
and with the seeding of straight line segments at points of locally
strongest attraction to obtain a dynamically consistent and numerically
robust approach to compute LCS. Extensions of these ideas to higher
dimensions are possible and will be communicated elsewhere.

\bibliographystyle{plain}

\end{document}